\documentclass[a4paper,twoside,11pt,reqno]{amsart}
\usepackage{amssymb}
\usepackage{amsmath}
\usepackage{amsfonts}

\newtheorem{theorem}{Theorem}

\newtheorem{definition}{Definition}
\newtheorem{example}{Example}

\newtheorem{proposition}{Proposition}
\newtheorem{remark}{Remark}

\title{   On Banachic Kernels and Approximation Theory}
\author{Marc Atteia} 
\email{marcatteia@orange.fr } 
\date{Version of \today}
\begin{document}
\maketitle 
\begin{abstract} 
In this paper, I generalize a previous one about hilbertian kernels and approximation theory
\end{abstract}
 
\section{ Definition of a banachic kernel.}

Let $\mathcal{E}$ be a locally convex vectorial space (lcs), $\mathcal{E}%
^{^{\prime }}$ its topological dual and $\left\langle \cdot ,\cdot
\right\rangle $ their duality bracket.
\subsection{ Hilbertian kernel.}
Let us suppose that $\mathcal{H}$ is a subspace of $\mathcal{E}$
and $\left( \mathcal{H},\left\langle \cdot \mid \cdot \right\rangle
\right) $ a hilbertian subspace of $\mathcal{E}$ .
\\
Then : \ $\mathcal{H}\subset \mathcal{E}$ and the injection $j$ from $%
\left( \mathcal{H},\left\langle \cdot \mid \cdot \right\rangle \right) $
into $\mathcal{E}$  is continuous.
\\
Let $\Lambda$ the duality mapping from $\mathcal{H}$ into $%
\mathcal{H}^{\prime }$ .
\\
The mapping $H=j\Lambda ^{t}j$ is called the hilbertian kernel of $%
\left( \mathcal{H},\left\langle \cdot \mid \cdot \right\rangle \right) $ \
relatively to $\mathcal{E}$ $.$
\\
So, we set : $\widetilde{H}=\Lambda ^{t}j$ .
\\
 One can prove that :
\begin{itemize}
\item[$\left( \ast \right)$]  $H$ is linear.
\item[$\left( \ast\ast \right)$]  $\widetilde{H}$$\mathcal{E}%
^{^{\prime }}$ is dense in $\left( \mathcal{H},\left\langle \cdot \mid
\cdot \right\rangle \right) $ .
\end{itemize}
\subsection{(Banachic) kernel of a reflexive and strictly
convex Banach space.}
 Let $\left( \mathcal{B},\left\Vert \cdot \right\Vert \right) $ a
reflexive and strictly convex Banach space, $\left( \mathcal{B}^{\prime
},\left\Vert \cdot \right\Vert _{\ast }\right) $  and 
$ \left( \mathcal{B}^{\prime \prime },\left\Vert \cdot \right\Vert
_{\ast \ast }\right) $  respectively,
its topological dual
and its toplogical bidual.
\\
 We suppose that $\mathcal{B}$ is a vectorial subspace of  $\mathcal{E}$ and that
the injection $j$ from $\left( \mathcal{B},\left\Vert \cdot
\right\Vert \right) $ into $\mathcal{E}$ is continuous.
\\
 So, we denote by :
\begin{itemize}
\item[$\left( \ast \right)$]  $\chi $ the canonical injection
from $\left( \mathcal{B},\left\Vert \cdot \right\Vert \right) $ into $%
\left( \mathcal{B}^{\prime \prime },\left\Vert \cdot \right\Vert _{\ast \ast
}\right) $ ;
\item[$\left( \ast\ast \right)$]
 $\psi $ an increasing function
from $ \mathbb{R}_{+}$ into $\mathbb{R} _{+}$ such that :
$$\psi \left( 0\right) =0\quad \textrm{and} \quad 
\lim_{\rho \rightarrow +\infty }\psi \left( \rho \right) =+\infty ;$$
\item[$\left( \ast\ast \ast\right)$]
$J$ $\left( \text{resp. }%
J_{\psi ^{-1}}^{\ast }\right) $ the duality mapping from  $\left( \mathcal{B},\left\Vert \cdot \right\Vert \right) $  into  $2^{\mathcal{B}^{\prime
}}$
$\left( \text{resp. from }\mathcal{B}^{\prime }\text{ into }2^{%
\mathcal{B}^{\prime \prime }\text{ }}\right) .$
\end{itemize}
 Thus :
\begin{equation*}
\left( \forall x\in \mathcal{B}\text{ , }J_{\psi }\left( x\right) =\left\{ 
\begin{array}{c}
x^{\ast }\in \mathcal{B}^{\prime }\text{ ; }\left\langle x^{\ast
},x\right\rangle =\left\Vert x^{\ast }\right\Vert _{\ast }\left\Vert
x\right\Vert \
\text{and : }\left\Vert x^{\ast }\right\Vert _{\ast }=\psi \left( \left\Vert
x\right\Vert \right)%
\end{array}%
\right\} \right)
\end{equation*}
Moreover we denote by  $Ban$ $\left( \mathcal{E}\right) $, the
set of all Banach spaces which are vectorial subspaces of $\mathcal{E}$
with a continuous injection into $\mathcal{E}$ .
\
The kernel of $\left( \mathcal{B},\left\Vert \cdot
\right\Vert \right) $ is the following mapping :%
\begin{equation*}
B=j\circ \chi ^{-1}\circ J_{\psi ^{-1}}^{\ast }\circ \text{ }^{t}j
\end{equation*}
  We set :  
\begin{equation*}
\widetilde{B}=\chi ^{-1}\circ J_{\psi ^{-1}}^{\ast }\circ \text{ }^{t}j
\end{equation*}
 Then :
\begin{itemize}
\item[$\left( \ast\right)$] 
 $B$  is a (multi-)mapping from  $
\mathcal{E}^{\prime }$ in $\mathcal{E}$ , non-linear generally.
\item[$\left( \ast\ast \right)$]  $\widetilde{B}$$\mathcal{E}%
^{\prime }$ is dense in $\left( \mathcal{B},\left\Vert \cdot \right\Vert
\right) $ .
\item[$\left( \ast\ast\ast \right)$] 
\begin{equation*}
\forall \text{ }x\in \mathcal{B}\text{ , }\forall e^{\prime }\in \mathcal{E}%
^{\prime }\text{ , }\left\langle jx,e^{\prime }\right\rangle =\psi \left(
\left\Vert \widetilde{B}e^{\prime }\right\Vert \right) .\left( \frac{d}{%
d\lambda }\left( \left\Vert \widetilde{B}e^{\prime }+\lambda x\right\Vert
\right) \right) _{\lambda =0}
\end{equation*}
\end{itemize}
\subsection{ The banachic kernel of an inf-compact potential.}
	\begin{definition}
	 A \textbf{potential} on $\mathcal{E}$ is a convex, even
functionnal vanishing at zero.
\
 A potential $\Phi $ on $\mathcal{E}$ is said $\sigma -\inf -$%
compact if there exists  $\lambda \in 
\mathbb{R}_{+}^{\ast }$ such that :
\begin{equation*}
S_{\lambda }\left( \Phi \right) =\left\{ e\in \mathcal{E}\text{ ; }\Phi
\left( e\right) \leq \lambda \right\} \text{ is compact in }\left( 
\mathcal{E},\sigma \left( \mathcal{E},\mathcal{E}^{\prime }\right)
\right)
\end{equation*}
\end{definition}
 A banachic kernel $B_{\Phi }$ of a $\sigma -\inf -$ compact
potential $\Phi $ is the subdifferential of
  the Legendre-Fenchel
 transform $\Phi ^{\ast }$ of $\Phi $ .
\\
So : $B_{\Phi }=\partial \Phi ^{\ast }.$
\begin{example}\end{example} 
 A hilbertian kernel is a Banachic kernel of a convenient
potential.
\begin{example}\end{example} 
\begin{equation*}
\forall e\in \mathcal{E},\text{ }\Phi \left( e\right) =\left\{ 
\begin{array}{c}
\omega \left( \left\Vert x\right\Vert \right) \text{ if }e=jx\text{ , }%
x\in \mathcal{B} \
+\infty \text{ otherwise} %
\end{array}%
\right\}
\end{equation*}  
where $\omega $  is a map from $ \mathbb{R}_{+}$ into $\mathbb{R}_{+}$ 
 derivable, even and such that : $\omega $ $\left( 0\right) =0$ .
\begin{proposition} 
 Let  $\Phi $  a potential on $\mathcal{E}$ and $E_{\Phi
} $ the vectorial subspace of $\mathcal{E}$ generated by
 $dom\Phi =\left\{ e\in \mathcal{E}\text{ ; }\Phi
\left( e\right) <+\infty \right\} $ .
  We set :%
\begin{equation*}
\forall x\in E_{\Phi }\text{ , }p_{\text{ }\Phi }\left( x\right) =\inf
\left\{ \lambda \in 
\mathbb{R}
_{+}^{\ast }\text{ : }\Phi \left( \frac{x}{\lambda }\right) \leq 1\text{}%
\right\}
\end{equation*}
 Then $p_{\text{ }\Phi }$ is a (semi-)norm on $%
E_{\Phi }$ .
\end{proposition}
\begin{definition}
   We say that $\Phi $ is reflexive if $\left(
E_{\Phi }\text{ , }p_{\text{ }\Phi }\right) $ is a reflexive normed space.
\end{definition}
\begin{proposition}
   $\Phi $ is reflexive iff the (convex) l.c.s.
regularized of $\Phi $ 
   is $\sigma -\inf -$compact on $\mathcal{E}$ .
\end{proposition}
\section{ First application}
 Let 
\begin{equation*}
a\in 
\mathbb{R}
_{+}^{\ast },m\in 
\mathbb{N}
^{\ast },p\in \left] 1,\infty \right[ \text{ and }p^{\ast \text{ }}%
\text{such that}:\frac{1}{p}+\frac{1}{p^{\ast }}=1
\end{equation*}
 We suppose that :%
\begin{eqnarray*}
\mathcal{E} &\mathcal{=}&%
\mathbb{R}
^{\left( -a,a\right) }\text{ . So :} \\
\mathcal{E}^{\prime } &=&\left\{ \sum\limits_{j\in I}\gamma _{j}\text{ }%
\delta _{t_{j}}\text{ ; }\gamma _{j}\in 
\mathbb{R}
\text{ , }t_{j}\in \left( -a,a\right) \text{, }Card\left( I\right)
<+\infty \right\}
\end{eqnarray*}
\subsection{   The space $\mathcal{A}_{p}$}
 Let us consider the following space :%
\begin{equation*}
\mathcal{A}_{p}=\left\{ x\in \mathcal{E}\text{ ; }\int_{-a}^{a}\left\vert
x^{\left( j\right) }\left( t\right) \right\vert ^{p}\text{ }dt<+\infty \text{
, }j=0,1,...m\right\}
\end{equation*}

 Then :%
\begin{equation*}
\mathcal{C}^{m}\left( -a,a\right) \subset \mathcal{A}_{p}\subset \mathcal{C}%
^{m-1}\left( -a,a\right)
\end{equation*}

 Given $e_{0}^{\prime },....,e_{m-1}^{\prime }\in \mathcal{E}%
^{\prime }$ , 
\begin{equation*}
\forall k\in \left\{ 0,1,...,\left( m-1\right) \right\} \text{ , }\forall
x\in \mathcal{A}_{p}\text{ , }\lambda _{k}\left( x\right) =\left\langle
jx,e_{k}^{\prime }\right\rangle
\end{equation*}
 Now, we denote by $\mathcal{P}_{m-1}\left( -a,a\right) $ , the
space of polynomials on $\left( -a,a\right) $ ,
 with degree less or equal to $\left( m-1\right) $,
 and we suppose that the restrictions of $\lambda _{0},\lambda
_{1},...,\lambda _{\left( m-1\right) }$ on $\mathcal{P}_{m-1}\left(
-a,a\right) $ are linearly independent.
 Then we set :%
\begin{equation*}
\forall x\in \mathcal{A}_{p}\text{ , }\left\Vert x\right\Vert _{\mathcal{A}%
_{p}}=\left( \sum_{k=0}^{m-1}\text{ }\left\vert \lambda _{k}\left( x\right)
\right\vert ^{p}+\int_{-a}^{a}\left\vert x^{\left( m\right) }\left( t\right)
\right\vert ^{p}\text{ }dt\right) ^{\frac{1}{p}}
\end{equation*}
 and :%
\begin{equation*}
\forall x\in \mathcal{A}_{p}\text{ , }\left\Vert x\right\Vert _{p}=\left(
\sum_{j=0}^{m}\left( \int_{-a}^{a}\left\vert x^{\left( j\right) }\left(
t\right) \right\vert ^{p}\text{ }dt\right) \right) ^{\frac{1}{p}}
\end{equation*}
 Then :
    $\left\Vert \cdot \right\Vert _{p}$ is a norm
on $\mathcal{A}_{p}$ and $\left( \mathcal{A}_{p},\left\Vert \cdot
\right\Vert _{p}\right) $ is a Banach space$.$
\begin{proposition}
If the functional $\lambda _{0},\lambda _{1},...,\lambda _{\left(
m-1\right) }$ are (linear and) continuous on $\left( \mathcal{A}%
_{p},\left\Vert \cdot \right\Vert _{p}\right) $
\begin{itemize}
  \item[(i)] The two norms $\left\Vert \cdot \right\Vert _{p}$ and
$\left\Vert \cdot \right\Vert _{\mathcal{A}_{p}}$ are equivalent
   (from the classical Banach theorem)
\\
 So, $\left( \mathcal{A}_{p},\left\Vert \cdot \right\Vert _{\mathcal{A}_{p}}\right) $ is a Banach space.
\item[(ii)] There exists $m$ polynomials $
P_{0},P_{1},...,P_{\left( m-1\right) }\in \mathcal{P}_{m-1}\left(-a,a\right) $
 such that :
\begin{equation*}
\lambda _{k}\left( P_{l}\right) =\delta _{kl}\text{ , }k,l\text{}\in
\left\{ 0,1...,\left( m-1\right) \right\}
\end{equation*}
\end{itemize}
  Moreover, we know that :
\begin{equation*}
\forall x\in \mathcal{C}^{m}\left( -a,a\right) \text{ , }\forall s\in
\left( -a,a\right) \text{ , }x\left( s\right) =\sum_{j=0}^{m-1}\frac{s^{j}%
}{j!}+\left( \int_{0}^{s}\frac{\left( s-t\right) _{+}^{m-1}}{\left(
m-1\right) !}.x^{\left( m\right) }\left( t\right) \text{ }dt\right)
\end{equation*}
\end{proposition}
 Thus, we can prove that :%
\begin{equation*}
\forall x\in \mathcal{C}^{m}\left( -a,a\right) \text{ , }\forall s\in
\left( -a,a\right) \text{ , }x\left( s\right) =\sum_{k=0}^{m-1}\lambda
_{k}\left( x\right) .P_{k}\left( s\right) +\left( \int_{-a}^{a}\Lambda
_{m}\left( s,t\right) .x^{\left( m\right) }\left( t\right) \text{ }dt\right)
\end{equation*}
 As $\mathcal{C}^{m}\left( -a,a\right) $ is dense in $\left( 
\mathcal{A}_{p},\left\Vert \cdot \right\Vert _{\mathcal{A}_{p}}\right) $ , \
the previous formula is true for any

$x\in \mathcal{A}_{p}$ , and the following proposition is true :
  \begin{proposition}
\begin{equation*}
\left( \mathcal{A}_{p},\left\Vert \cdot \right\Vert _{\mathcal{A}%
_{p}}\right) \in Ban\left( \mathcal{E}\right)
\end{equation*}
\end{proposition}
 Let us set : 
\begin{equation*}
\mathcal{R}_{p}=\mathcal{P}_{m-1}\left( -a,a\right) \text{ and }%
\forall P\in \mathcal{R}_{p}\text{ , }\left\Vert P\right\Vert _{\mathcal{R}%
_{p}}=\left( \sum_{k=0}^{m-1}\text{ }\left\vert \lambda _{k}\left( P\right)
\right\vert ^{p}\right) ^{\frac{1}{p}}
\end{equation*}
Then  $\left\Vert \cdot \right\Vert _{\mathcal{R}_{p}}$ is a
norm on $\mathcal{R}_{p}$ .
\\
So, $\mathcal{A}_{p}$ is the algebra\"{\i}c direct sum of $%
\mathcal{R}_{p}$ and :%
\begin{equation*}
\mathcal{C}_{p}=\left\{ x\in \mathcal{E}\text{ ; }\int_{-a}^{a}\left\vert
x^{\left( m\right) }\left( t\right) \right\vert ^{p}\text{ }dt<+\infty \text{
with }\lambda _{k}\left( x\right) =0\text{ , }k\in \left\{
0,1...,\left( m-1\right) \right\} \right\}
\end{equation*}
Let us set :%
\begin{equation*}
\forall x\in \mathcal{C}_{p}\text{ , }\left\Vert x\right\Vert _{\mathcal{%
C}_{p}}=\left( \int_{-a}^{a}\left\vert x^{\left( m\right) }\left( t\right)
\right\vert ^{p}\text{ }dt\right) ^{\frac{1}{p}}
\end{equation*} 
Then:
\begin{equation*}
\left( \mathcal{C}_{p},\left\Vert \cdot \right\Vert _{\mathcal{C}%
_{p}}\right) \text{ is a Banach space}
\end{equation*}
As $dim\left( \mathcal{R}_{p}\right) $ is finite, 
\begin{equation*}
\left( \mathcal{A}_{p},\left\Vert \cdot \right\Vert _{\mathcal{A}%
_{p}}\right) \text{ is the direct topological sum of }\left( \mathcal{R}%
_{p},\left\Vert \cdot \right\Vert _{\mathcal{R}_{p}}\right) \text{ and }%
\left( \mathcal{C}_{p},\left\Vert \cdot \right\Vert _{\mathcal{C}_{p}}\right)
\end{equation*}
 \subsection{  The banachic kernel of $\left( \mathcal{%
A}_{p},\left\Vert \cdot \right\Vert _{\mathcal{A}_{p}}\right) $ }

 $\left( \ast \right) $ Given $\mathcal{B}\in Ban$ $\left( \mathcal{E%
}\right) ,$ we set :%
\begin{equation*}
\forall x\in \mathcal{E}\text{ , }f_{\mathcal{B}}\left( e\right) =\left\{ 
\begin{array}{c}
\omega _{\mathcal{B}}\left( \left\Vert \cdot \right\Vert _{\mathcal{B}%
}\right) \text{ if }e=jx,\text{ }x\in \mathcal{B} \
+\infty \text{ otherwise}%
\end{array}%
\right\}
\end{equation*}

 In the following, $\mathcal{B}$ will be identified to $\mathcal{A}%
_{p}$ , $\mathcal{R}_{p}$ and $\mathcal{C}_{p}$ and we suppose that :%
\begin{equation*}
\omega _{\mathcal{B}}\left( \rho \right) =\rho ^{p}\text{ when }\rho \in 
\mathbb{R}^*
\end{equation*}

 Let us set :%
\begin{equation*}
\forall \rho \in \mathbb{R}^*\text{ , }\alpha _{p}\left( \rho \right)
=\left\vert \rho \right\vert ^{p-1}.signe\left( \rho \right) \text{ }%
\implies \left( \alpha _{p}\right) ^{-1}=\alpha _{p^{\ast }}
\end{equation*}

 Then : 
\begin{eqnarray*}
\forall x,y &\in &\mathcal{A}_{p}\text{ , }\forall \mu \in \mathbb{R}\text{
, }\left[ \frac{\partial }{\partial \mu }\left( f\left( y+\mu x\right)
\right) \right] _{\mu =0} \\
&=&\sum_{k=0}^{m-1}\text{ }\alpha _{p}.\lambda _{k}\left( y\right) .\lambda
_{k}\left( x\right) +\int_{-a}^{a}\alpha _{p}.y^{\left( m\right) }\left(
t\right) .x^{\left( m\right) }\left( t\right) \text{ }dt
\end{eqnarray*}

 Let $A_{p}$ (resp. $R_{p}$ and $C_{p}$) the banachic kernel
of $\left( \mathcal{A}_{p},\left\Vert \cdot \right\Vert _{\mathcal{A}%
_{p}}\right) $

 (resp. $\left( \mathcal{R}_{p},\left\Vert \cdot \right\Vert _{%
\mathcal{R}_{p}}\right) $ , $\left( \mathcal{C}_{p},\left\Vert \cdot
\right\Vert _{\mathcal{C}_{p}}\right) $) relatively to $f_{\mathcal{A}%
_{p}}$ (resp. $f_{\mathcal{R}_{p}}$ , $f_{\mathcal{C}_{p}}$).

  But, $\forall x\in $ $\mathcal{A}_{p}$ , $x=x_{1}+x_{2}$ , 
$x_{1}\in \mathcal{R}_{p}$ , $x_{2}$ $\in \mathcal{C}_{p}$ (this
decomposition is unique).

  We deduce that : 
\begin{equation*}
\forall x\in \mathcal{A}_{p}\text{ , }f_{_{\mathcal{A}_{p}}}\left(
x\right) =f_{_{\mathcal{R}_{p}}}\left( x\right) +f_{_{\mathcal{C}%
_{p}}}\left( x\right) \text{}
\end{equation*}
\begin{proposition}
\begin{itemize}
   \item [(i)] 
\begin{eqnarray*}
\forall e^{\prime },\text{ }f^{\prime } &\in &\mathcal{E}^{\prime },\text{ }%
\left\langle \mathcal{R}_{p}e^{\prime },f^{\prime }\right\rangle
=\sum_{k=0}^{m-1}\text{ }\alpha _{p^{\ast }}\left\langle jP_{k},e^{\prime
}\right\rangle .\left\langle jP_{k},f^{\prime }\right\rangle \\
\left\langle \mathcal{C}_{p}e^{\prime },f^{\prime }\right\rangle
&=&\int_{-a}^{a}\left\langle j\Lambda _{m}\left( \cdot ,\theta \right)
,f^{\prime }\right\rangle .\text{ }\alpha _{p^{\ast }}\left\langle j\Lambda
_{m}\left( \cdot ,\theta \right) ,e^{\prime }\right\rangle \text{ }d\theta
\end{eqnarray*}
\item[(ii)] $\displaystyle A_{p}=R_{p} +C_{p}$ .
\end{itemize}
    So :
\begin{equation*}
\forall s,t\text{ }\in \left( -a,a\right) \text{ }pp\text{ , }%
\left\langle A_{p}\delta _{s},\delta _{t}\right\rangle =\sum_{k=0}^{m-1}%
\text{ }P_{k}\left( t\right) .\alpha _{p^{\ast }}P_{k}\left( s\right)
+\int_{-a}^{a}\Lambda _{m}\left( t,\theta \right) \text{ }\alpha _{p^{\ast
}}\Lambda _{m}\left( s,\theta \right) \text{ }d\theta
\end{equation*}
    \end{proposition}
    \begin{remark}           \end{remark} 
   One can verify that : $A_{2},$ $R_{2}$ $,C_{2}$ \
are hilbertian kernels.
\subsection{ Some properties of $C_{p}$ }
\quad \\\  (i) Let $a\in R_{+}^{\ast }$ $,$ $p,q\in \left] 1,\infty %
\right[ $ with $p<q$ ; then : $\mathcal{C}_{q}\subset \mathcal{C}_{p}$ \
and 

$  \mathcal{C}_{q}$ is dense in $\mathcal{C}_{p}$ .

  Now, we set : 
\begin{equation*}
\forall s,t\text{ }\in \left( -a,a\right) \text{ }pp\text{ , }C_{p}\left(
t,s\right) =\left\langle C_{p}\delta _{s},\delta _{t}\right\rangle
\end{equation*}

  From previous results we deduce that :%
\begin{equation*}
\forall s,t\text{ }\in \left( -a,a\right) \text{ }pp\text{ , }\frac{%
\partial ^{m}}{\partial t^{m}}\left( \mathcal{C}_{q}\left( t,s\right)
\right) =\alpha _{1+\frac{p-1}{q-1}}\left( \frac{\partial ^{m}}{\partial
t^{m}}\left( \mathcal{C}_{p}\left( t,s\right) \right) \right)
\end{equation*}

   This formula is symmetric relatively to $p$ and $%
q $ because :%
\begin{equation*}
\left( \alpha _{1+\frac{p-1}{q-1}}\right) ^{-1}=\alpha _{1+\frac{q-1}{p-1}}
\end{equation*}

   Thus, this formula is true for any $p,q\in \left]
1,\infty \right[ $ .

   So, when $p=2$ , we have :%
\begin{eqnarray*}
\forall q &\in &\left] 1,\infty \right[ \text{ , }\forall s,t\text{ }\in
\left( -a,a\right) \text{ }pp\text{ , }\frac{\partial ^{m}}{\partial t^{m}}%
\left( \mathcal{C}_{q}\left( t,s\right) \right) =\alpha _{q^{\ast }}\left( 
\frac{\partial ^{m}}{\partial t^{m}}\left( \mathcal{C}_{2}\left( t,s\right)
\right) \right) \text{ , } \\
\frac{1}{q}+\frac{1}{q^{\ast }} &=&1\text{}
\end{eqnarray*}

  As, $\alpha _{p}+\alpha _{p^{\ast }}=1$ and as $\mathcal{C}%
_{2}$ is a hilbertian kernel,

  therefore a symmetric kernel, we have :%
\begin{eqnarray*}
\text{ }\forall s,\theta \text{ } &\in &\left( -a,a\right) \text{ }pp\text{
,} \\
\mathcal{C}_{p}\left( \theta ,s\right) &=&\int_{-a}^{a}\left( \frac{\partial
^{m}}{\partial t^{m}}\left( \mathcal{C}_{2}\left( t,s\right) \right) \right) 
\text{ }\alpha _{p^{\ast }}\left( \frac{\partial ^{m}}{\partial t^{m}}\left( 
\mathcal{C}_{2}\left( \theta ,t\right) \right) \right) \text{ }dt \\
\text{with } &\text{: }&\mathcal{C}_{p}\left( s,s\right)
=\int_{-a}^{a}\left( \frac{\partial ^{m}}{\partial t^{m}}\left( \mathcal{C}%
_{2}\left( t,s\right) \right) \right) ^{p^{\ast }}\text{ }dt
\end{eqnarray*}

  (ii) When $\mathcal{E}=\mathcal{D}^{\prime }\left(
-a,a\right) $ and $\mathcal{E}^{\prime }=\mathcal{D}\left( -a,a\right) $ ,
we have :%
\begin{eqnarray*}
\forall t\text{ } &\in &\left( -a,a\right) \text{ }pp\text{ , }\forall
\varphi ,\psi \in \mathcal{E}\text{ ,} \\
\left( \ast \right) \text{ }\left( -1\right) ^{m}\alpha _{p}\left( \frac{%
\partial ^{m}}{\partial t^{m}}\left( \widetilde{\mathcal{C}_{p}}\varphi
\right) \left( t\right) \right) &=&\varphi \left( t\right) \text{ with : }
\\
\lambda _{k}\left( \widetilde{\mathcal{C}_{p}}\varphi \right) &=&0\text{ , 
}k=0,1,...,\left( m-1\right)
\end{eqnarray*}

\begin{eqnarray*}
\left( \ast \ast \right) \text{ }\left\langle \mathcal{C}_{p}\psi ,\varphi
\right\rangle &=&\int_{-a}^{a}\left( \widetilde{\mathcal{C}_{p}}\psi \right)
\left( t\right) \varphi \left( t\right) \text{ }dt \\
&=&\int_{-a}^{a}\left( \frac{\partial ^{m}}{\partial t^{m}}\left( \widetilde{%
\mathcal{C}_{2}}\varphi \right) \left( t\right) \right) \text{ }\alpha
_{p^{\ast }}\left( \frac{\partial ^{m}}{\partial t^{m}}\left( \widetilde{%
\mathcal{C}_{2}}\psi \right) \left( t\right) \right) \text{ }dt
\end{eqnarray*}

    and :%
\begin{equation*}
\left\langle \mathcal{C}_{p}\varphi ,\varphi \right\rangle
=\int_{-a}^{a}\left( \frac{\partial ^{m}}{\partial t^{m}}\left( \widetilde{%
\mathcal{C}_{2}}\varphi \right) \left( t\right) \right) ^{p^{\ast }}\text{ }%
dt
\end{equation*}
 \subsection{ Banachic kernel and the Sard's factorization theorem%
}
  Let $D^{m}$ the derivative of order $m$ ; that is a
linear and continuous mapping
 from $\left( \mathcal{A}_{p},\left\Vert \cdot \right\Vert _{\mathcal{A%
}_{p}}\right) $ onto $L^{p}\left( a,b\right) $ .
\\
  Given $e^{\prime }\in \mathcal{E}^{\prime }$ , we shall
denote $v_{e^{\prime }}$ the mapping from $\mathcal{A}_{p}$ into $%
\mathbb{R}$ such that :%
\begin{equation*}
\forall x\in \mathcal{A}_{p}\text{ , }v_{e^{\prime }}\left( x\right)
=\left\langle j\left( x-\sum_{k=0}^{m-1}\text{ }\lambda _{k}\left( x\right)
.P_{k}\right) ,e^{\prime }\right\rangle
\end{equation*}
 So, we have the following scheme :
\\          
$$  \begin{array}{ccc}
&\mathcal{A}_{p}& \\
   &\swarrow    \, \searrow &       \\
L^{p}\left( a,b\right) &\longrightarrow &\mathbb{R}
\end{array}
$$ 
From the Sard's factorization theorem we deduce that there
exists $G\in L^{p^{\ast }}\left( a,b\right) $ such that :
\begin{equation*}
\text{}v_{e^{\prime }}\left( x\right) =\int_{a}^{b}G\left( t\right)
.x^{\left( m\right) }\left( t\right) \text{ }dt
\end{equation*}
 thus :
\begin{equation*}
G=\alpha _{p}\left( \left( \widetilde{\mathcal{C}_{p}}e^{\prime }\right)
^{\left( m\right) }\right)
\end{equation*}
 \section{ Banachic B-splines}
\subsection{About an approximation problem in a Banach
space}
  Let $n\in \mathbb{N}$ , $e_{0}^{\prime },e_{1}^{\prime
},...,e_{n}^{\prime }$ , $\alpha _{0},\alpha _{1},...,\alpha _{n}$ $%
\in \mathbb{R}$ and : 
\begin{equation*}
\Gamma =\left\{ e\in \mathcal{E}\text{ ; }\left\langle e,e_{k}^{\prime
}\right\rangle =\alpha _{k}\text{ , }k=0,1,...,n\right\}
\end{equation*}
   Then $\Gamma $ is an hyperplane (with finite
codimension) which is closed
  in $\left( \mathcal{E},\sigma \left( \mathcal{E}\text{ },%
\mathcal{E}\text{ }^{\prime }\right) \right) $.
\\
  Let $\Phi $ be a $\sigma -\inf -$ compact potential
on $\mathcal{E}$ .
   Let us consider the following problem :%
\begin{equation*}
\Pi \text{ : }\inf \left\{ \Phi \left( e\right) \text{ ; }e\in \Gamma
\right\}
\end{equation*}
\begin{proposition}
  Let us suppose that :%
\begin{equation*}
\inf \left\{ \Phi \left( e\right) \text{ ; }e\in \Gamma \right\} =\mu \in 
\mathbb{R}
\end{equation*}
  Then $\Pi $ has a solution (which is unique when $\Phi $ \
is strictly convex).
\end{proposition}
\subsection{ Characterization of a solution $\overline{e}$ of $\Pi $ }
\begin{theorem} 
  Let $\Gamma _{0}\left( \mathcal{E}\right) $ be the set of
proper convex l.c.s. functional on $\mathcal{E}$ ,
  which are not identical at +$\infty $ ,
 $ \Delta $ a convex subset of $\mathcal{E}$ and $u\in
\Gamma _{0}\left( \mathcal{E}\right) .$
\\
  Let us consider the following problem :%
\begin{equation*}
\Pi \text{ : }\inf \left\{ u\left( e\right) \text{ ; }e\in \Delta \right\}
\end{equation*}
  If the following hypotheses are verified :
\begin{itemize}
 \item[$\left( \ast \right) $]
  $\inf \left\{ u\left( e\right) 
\text{ ; }e\in \Delta \right\} =\mu \in \mathbb{R}$
\item[$\left( \ast \ast \right)$] $u$ is finite and
continuous at a point of $\Delta $
\\
     or $u$ is finite at an interior point
of $\Delta $
\end{itemize}
  then :
\begin{itemize}
\item[(i)]
\begin{equation*}
{\mu}%
=\max \left\{ -u^{\ast }\left( -e^{\prime }\right) -h\left( e^{\prime };%
\overline{\Delta }\right) \text{ };\text{ }e^{\prime }\in dom\left( h\left(
\cdot \text{ ; }\overline{\Delta }\right) \right) \right\}
\end{equation*}
   where $u^{\ast }$ is the dual of $u$ , $%
\overline{\Delta }$ is the closure of $\Delta $~in $\left( \mathcal{E%
},\sigma \left( \mathcal{E}\text{ },\mathcal{E}\text{ }^{\prime }\right)
\right) $

   and $h\left( \cdot \text{ };\overline{\Delta }\right)
=\left( \delta \left( \cdot \text{ };\overline{\Delta }\right) \right)
^{\ast }$
\item[(ii)]
 Moreover, if 
\begin{equation*}
\left\{ e\in \mathcal{E}\text{ ; }u\left( e\right) +\delta \left( e\text{ }%
;\Delta \right) \leq 
{\mu}%
\right\}
\end{equation*}
   is a non-void set
\\
   or, if%
\begin{equation*}
\exists \text{ }e_{0}^{\prime }\in \mathcal{E}\text{ }^{\prime }\text{ ; }%
\partial u^{\ast }\left( -e_{0}^{\prime }\right) \cap \partial h\left(
e_{0}^{\prime }\text{ };\overline{\Delta }\right)
\end{equation*}
   is a non-void set,
\\
   then :%
\begin{eqnarray*}
&&\left\{ e\in \mathcal{E}\text{ ; }u\left( e\right) +\delta \left( e\text{
};\Delta \right) \leq 
{\mu}%
\right\} \\
&=&\cup \left( \left( \text{ }\partial u^{\ast }\left( -e^{\prime }\right)
\cap \partial h\left( e^{\prime }\text{ };\overline{\Delta }\right) \right) 
\text{ };\text{ }e^{\prime }\in \mathcal{E}\text{ }^{\prime }\text{ ;}\right)
\end{eqnarray*}
\end{itemize}
\end{theorem}
\begin{proposition}
  Let us suppose that $\Phi $ is finite and continuous at a
point in $\Gamma .$
\\
  Let $\overline{e}$ a solution of $\Pi $ . Then : 
\begin{equation*}
\overline{e}\in B_{\Phi }\left( \sum_{j=0}^{n}\overline{\lambda _{j}}\text{ }%
e_{j}^{\prime }\right) \text{ where }\overline{\lambda }_{0},\overline{%
\lambda }_{1},...,\overline{\lambda }_{n}\in \mathbb{R}\text{ and }%
\overline{e}\in \Delta \text{ and }B_{\Phi }=\partial \Phi ^{\ast }
\end{equation*}
 \end{proposition}
\subsection{  An application}
  I shall use below the same hypotheses and notations as in the
previous paragraphs.
\\
  Given $e_{0}^{\prime },...,e_{m-1}^{\prime },e_{m}^{\prime
}...,e_{n}^{\prime }\in \mathcal{E}^{\prime }$ and $%
a_{0},...,a_{m},a_{m-1,}...,a_{n}\in \mathbb{R},$ we consider the
 following problem :%
\begin{equation*}
\Pi _{0}\text{ : }\inf \left\{ \text{ ;}\int_{-a}^{a}\left\vert x^{\left(
m\right) }\left( t\right) \right\vert ^{p}\text{ }dt\text{ ; }x\in \mathcal{A%
}_{p}\text{ },\text{ }\left\langle jx,e_{l}^{\prime }\right\rangle =\alpha
_{l}\text{ },\text{ }l\in \left\{ 0,...,n\right\} \right\}
\end{equation*}
$\Pi _{0}$ is equivalent to :%
\begin{equation*}
\Pi _{1}\text{ : }\inf \left\{ \text{ ;}\int_{-a}^{a}\left\vert x^{\left(
m\right) }\left( t\right) \right\vert ^{p}\text{ }dt\text{ }+\sum_{k=0}^{m-1}%
\text{ }\left\vert \lambda _{k}\left( x\right) \right\vert ^{p}\text{; }x\in 
\mathcal{A}_{p}\text{ },\text{ }\left\langle jx,e_{l}^{\prime
}\right\rangle =\alpha _{l}\text{ },\text{ }l\in \left\{ 0,...,n\right\}
\right\}
\end{equation*}
  Thus, we deduce that :%
\begin{equation*}
\sigma _{p}=\widetilde{A}\left( \sum_{l=0}^{n}\mu _{l}\text{ }e_{l}^{\prime
}\right) _{p}\text{ ; }\mu _{0},\mu _{1},...,\mu _{n}\in \mathbb{R}\text{
with }\left\langle j\sigma _{p},e_{l}^{\prime }\right\rangle =\alpha
_{l}\text{ },\text{ }l\in \left\{ 0,...,n\right\}
\end{equation*}
  So :%
\begin{eqnarray*}
\forall t\text{ } &\in &\left( -a,a\right) \text{ }pp\text{ ,} \\
\sigma _{p}\left( t\right) &=&\sum_{k=0}^{m-1}P_{k}\left( t\right) .\alpha
_{p^{\ast }}\left[ \sum_{l=0}^{n}\mu _{l}\text{}\left\langle
jP_{k},e_{l}^{\prime }\right\rangle \right] \\
&&+\int_{-a}^{a}\Lambda _{m}\left( t,\theta \right) \text{ }\alpha _{p^{\ast
}}\left[ \sum_{l=0}^{n}\mu _{l}\text{}\left\langle j\Lambda _{m}\left(
\cdot ,\theta \right) ,e_{l}^{\prime }\right\rangle \right] \text{ }d\theta
\\
\text{with }\left\langle j\sigma _{p},e_{l}^{\prime }\right\rangle
&=&\alpha _{l}\text{ },\text{ }l\in \left\{ 0,...,n\right\}
\end{eqnarray*}
 Then, to calculate $\mu _{0},\mu _{1},...,\mu _{n}$ we
must solve a non-linear system of equations.
 \\
\textbf{  Properties of $\protect\sigma _{p}$ :}
\\
   Let $\tau _{2}$ the solution of the following
problem :%
\begin{equation*}
\Pi _{2}\text{ : }\inf \left\{ \text{ ;}\int_{-a}^{a}\left\vert x^{\left(
m\right) }\left( t\right) \right\vert ^{2}\text{ }dt\text{ ; }x\in \mathcal{A%
}_{2}\text{ },\text{ }\left\langle jx,e_{l}^{\prime }\right\rangle =\beta
_{l}\text{ },\text{ }l\in \left\{ 0,...,n\right\} \right\}
\end{equation*}
  Then :%
\begin{eqnarray*}
\forall t\text{ } &\in &\left( -a,a\right) \text{ }pp\text{ ,} \\
\tau _{2}\left( t\right) &=&\sum_{k=0}^{n}\nu
_{k}\sum_{l=0}^{m-1}P_{l}\left( t\right) .\left\langle jP_{l},e_{k}^{\prime
}\right\rangle \\
&&+\int_{-a}^{a}\Lambda _{m}\left( t,\theta \right) \text{ }%
\sum_{l=0}^{n}\left\langle j\Lambda _{m}\left( \cdot ,\theta \right)
,e_{k}^{\prime }\right\rangle \text{ }d\theta \\
\text{with }\left\langle j\tau _{2},e_{l}^{\prime }\right\rangle
&=&\beta _{l}\text{ },\text{ }l\in \left\{ 0,...,n\right\}
\end{eqnarray*} 
  Let us suppose that we have choose $\beta _{0},...,\beta
_{n}\in \mathbb{R}$ , such that :

   $\nu _{l}=\mu _{l}$ , $l=0,1,...,n$ .
  Then :%
\begin{eqnarray*}
\forall t\text{ } &\in &\left( -a,a\right) \text{ }pp\text{ ,} \\
\tau _{2}^{\left( m\right) }\left( t\right) &=&\sum_{l=0}^{n}\mu _{l\text{ }%
}\left\langle j\Lambda _{m}\left( \cdot ,t\right) ,e_{l}^{\prime
}\right\rangle =\alpha _{p}\left[ \sigma ^{\left( m\right) }\left( t\right) %
\right]
\end{eqnarray*}
  Thus : 
\begin{equation*}
\sigma _{p}^{\left( m\right) }\left( t\right) =\alpha _{p^{\ast }}\left(
\tau _{2}^{\left( m\right) }\left( t\right) \right)
\end{equation*}
  So, if $e_{l}^{\prime }=\delta _{t_{l}}$ , $t_{l}$ $\in
\left( -a,a\right) $ , $l=0,1,...,n$ , we have :%
\begin{equation*}
\sigma _{p}^{\left( m\right) }\left( t\right) =\alpha _{p^{\ast }}\left(
\sum_{l=0}^{n}\mu _{l\text{ }}\Lambda _{m}\left( t_{l}\text{ },t\right)
\right)
\end{equation*}
\subsection{ Banachic B-splines}

  Let $L_{k,m-1}$ be the polynomial with degree less or
equal to $\left( m-1\right) $ such that : 
\begin{equation*}
L_{k,m-1}\left( t_{l}\right) =\delta _{k,l\text{ }}\text{ },\text{ }%
k,l=0,1,...,\left( m-1\right)
\end{equation*}

  Then :%
\begin{eqnarray*}
\forall s,t\text{ } &\in &\left( -a,a\right) \text{ }pp \\
\Lambda _{m}\left( s\text{ },t\right) &=&\frac{1}{\left( m-1\right) !}%
\left[ \left( s-t\right) _{+}^{m-1}-\sum_{k=0}^{n}\left( t_{k}-t\right)
_{+}^{m-1}L_{k,m-1}\left( s\right) \right]
\end{eqnarray*}

  Now, let us suppose that $n\geq m$.

  Let us choose $\mu _{0},...,\mu _{n}\in \mathbb{R}$ ,
such that :%
\begin{equation*}
\sum_{l=0}^{n}\mu _{l\text{ }}f\left( t_{l}\right) \text{ is equal to
the divided difference of order }m,\text{ at the point }\tau \in \left(
-a,a\right)
\end{equation*}

  So :%
\begin{eqnarray*}
\forall t\text{ } &\in &\left( -a,a\right) \text{ }pp\text{ ,} \\
\tau _{2}^{\left( m\right) }\left( t\right) &=&\sum_{l=0}^{n}\mu _{l}\frac{1%
}{\left( m-1\right) !}\left[ \left( t_{l}-t\right)
_{+}^{m-1}-\sum_{k=0}^{n}\left( t_{k}-t\right) _{+}^{m-1}L_{k,m-1}\left( \
t_{l}\right) \right] \\
&=&\sum_{l=0}^{n}\mu _{l}\frac{\left( t_{l}-t\right) _{+}^{m-1}}{\left(
m-1\right) !}
\end{eqnarray*}

  Thus :

   $\tau _{2}^{\left( m\right) }$ is a classical
polynomial B-spline of degree equal at $\left( m-1\right) $ .

  The B-spline $\tau _{2}^{\left( m\right) }$ is equal to
zero on the set $\left( -\infty ,a_{1}\right) \cup \left(
a_{2},+\infty \right) $ 
  with $a_{1},a_{2}\in \left( -a,a\right) $ and $\
a_{1}<a_{2}.$

  (We remark moreover that $\tau _{2}^{\left( m\right) }$ \
is an optimization spline function.)

As $\alpha _{p^{\ast }}=0$ , because $p^{\ast }>1$ ,

    $\sigma _{p}^{\left( m\right) }$ is equal to
zero on the set $\left( -\infty ,a_{1}\right) \cup \left(
a_{2},+\infty \right) .$
\begin{definition} 
   I say that $\sigma _{p}^{\left( m\right) }$ is a
banachic B-spline.
\end{definition}
  Let us suppose that : $t_{l+1}-t_{l}=h$ , $%
l=0,1,...,\left( n-1\right) $ ;
\\
  Let us set :%
\begin{equation*}
D_{h}^{m}f=\sum_{l=0}^{n}\mu _{l}\text{ }f\left( t_{l}\right) \text{ and 
}Q_{h}^{m,2}=\tau _{2}^{\left( m\right) }\text{ , }Q_{h}^{m,p}=\sigma
_{p}^{\left( m\right) }
\end{equation*}
  Then :%
\begin{equation*}
Q_{h}^{m,2}=\alpha _{p}\left( Q_{h}^{m,p}\right) \text{ with : }\alpha
_{p}\left( \rho \right) =\left\vert \rho \right\vert ^{p-1}sign\left( \rho
\right)
\end{equation*}
  From the classical properties of $Q_{h}^{m,2}$ , we deduce
easily those of $Q_{h}^{m,p}$ .
  So,
\begin{proposition} 
\quad\\
\begin{itemize}
\item[(i)] $0\leq Q_{h}^{m,p}\left( s,t\right) \leq h^{1-p}$
\item[(ii)] $\int_{\mathbb{R}}$ $\left( Q_{h}^{m,p}\left(
s,t\right) \right) ^{p-1}ds=1$ and $\sum_{j\in \mathbb{Z}}$ $\left(
Q_{h}^{m,p}\left( jh,t\right) \right) ^{p-1}ds=h^{-1}$
\item[(iii)] $\forall y\in \mathcal{C}^{m}\left( \mathbb{R}%
\right) $ , $\int_{\mathbb{R}}$ $y^{\left( m\right) }\left( s\right)
\left( Q_{h}^{m,p}\left( s,t\right) \right) ^{p-1}ds=D_{h}^{m}y\left(
t\right) $   
\end{itemize}
\end{proposition}
    Other properties can be found in the book of the author
:

  "Hilbertian kernels and spline functions".
\section{  Extensions}
\subsection{}
$\ast$ $M$ is a potential such that $M^{\ast }$ is derivable or

 when $M$ is convex, positive, such that $M\left( 0\right) =0$
and $M^{\ast }$ is derivable,

 we can consider the case where :%
\begin{equation*}
\frac{\partial ^{m}}{\partial t^{m}}\left( \widetilde{C}_{p}e^{\prime
}\right) \left( t\right) =\left( M^{\ast }\right) ^{\prime }\left(
\left\langle j\Lambda _{m}\left( \cdot ,t\right) ,e^{\prime }\right\rangle
\right)
\end{equation*}

 Then, what space is associated to $M$ ?
\subsection{  Calculus}
 of the banachic kernels of the spaces $\mathcal{B}_{p,k}$ when
$p\in \left] 1,+\infty \right[ $
where $\mathcal{B}_{p,k}$ is the space of distributions $%
x\in S^{\prime }\left( \mathbb{R}\right) $ such that :%
\begin{equation*}
\left\Vert x\right\Vert _{p,k}=\left[ \left( 2\pi \right) ^{-n}\int_{\mathbb{%
R}^{n}}\left\vert k\left( \zeta \right) .\widehat{x}\left( \zeta \right)
\right\vert ^{p}d\zeta \right] ^{\frac{1}{p}}
\end{equation*}
 $ k$ is a mapping from $\mathbb{R}^{n}$ into $\mathbb{R}_{+}$
such that :%
\begin{equation*}
\forall \zeta ,\eta \in \mathbb{R}^{n}\text{ , }k\left( \zeta +\eta
\right) \leq \left( 1+c\left\vert \zeta \right\vert ^{\nu }\right) .k\left(
\eta \right)
\end{equation*}
There are also many other extensions to banachic kernels of the
properties of hilbertian
kernels (wavelets, fractals, vectorial case, etc.)
\section{Banachic kernels and partial differntial equations.
An example.}

  Let $K$ a triangle in $\mathbb{R}^{2}$ which vertices
are $A=\left( 0,0\right) $ , $B=\left( 1,0\right) $ , $C=\left(
0,1\right) $ .

 We denote by $int\left( K\right) $ the interior of $K$ .

  Let $\mathcal{E}=\mathcal{D}^{\prime }\left( int\left( \
K\right) \right) $ and $\mathcal{E}^{\prime }=\mathcal{D}\left( \
int\left( K\right) \right) $ .

 Given $p\in \left] 1,+\infty \right[ $ , let us consider the
Banach space $\mathcal{B}$ of continuous functions on $K$
 which are equal to zero at $A,B,C\,$embedded with the
following norm :%
\begin{equation*}
v\rightarrow \sum_{\left\vert \alpha \right\vert =1}\left(
\int_{K}\left\vert D^{\alpha }v\left( t\right) \right\vert ^{p}dt\right) ^{%
\frac{1}{p}}
\end{equation*}

  Then , using the same notations as those which are defined at
the paragraph 1.2,

 we have :%
\begin{eqnarray*}
\forall \varphi &\in &\mathcal{E}^{\prime },\text{ }\forall v\in 
\mathcal{B}\text{ ,} \\
\left\langle jv,\varphi \right\rangle &=&\left\Vert \widetilde{B}\varphi
\right\Vert ^{p-1}.\left( \frac{d}{d\lambda }\left( \left\Vert \widetilde{B}%
\left( \varphi +\lambda v\right) \right\Vert \right) \right) _{\lambda =0}
\end{eqnarray*}

  But :%
\begin{equation*}
\frac{d}{d\lambda }\left( \left\Vert u+\lambda v\right\Vert \right)
=\left\Vert u\right\Vert ^{1-p}.\sum_{\left\vert \beta \right\vert =1}\left(
\int_{K}\left\vert D^{\beta }u\left( t\right) \right\vert ^{p-1}.sign\left(
D^{\beta }u\left( t\right) \right) .D^{\beta }v\left( t\right) dt\right)
\end{equation*}

 Now, we set :%
\begin{equation*}
\forall u\in \mathcal{E}\text{ , }Lu=\sum_{\left\vert \beta \right\vert
=1}\left( \int_{K}D^{\beta }\left( \left\vert D^{\beta }u\right\vert
^{p-1}.sign\left( D^{\beta }u\right) \right) .\right)
\end{equation*}

  So, the banachic kernel $B$ of $\mathcal{B}$ is the
solution of the following non linear differential

 problem :%
\begin{equation*}
\forall \varphi \in \mathcal{E}^{\prime },\text{ }L\left( B\varphi \right)
=\varphi \text{ with }\left( B\varphi \right) \left( A\right)
=\left( B\varphi \right) \left( B\right) =\left( B\varphi \right) \left(
C\right) =0
\end{equation*}%
\textbf{References. }
 \vspace{0.3cm}
 
 [1]  M. Atteia, J. Audounet. inf compact potentials and banachic kernels. Lectures notes in Mathematics, 991, pages 7--27, 1983.

\end{document}